\documentclass{amsart}
\usepackage{latexsym}
\usepackage{amssymb}
\usepackage{amsmath}

\begin{document}


\newtheorem{theorem}{Theorem}
\newtheorem{problem}{Problem}
\newtheorem{definition}{Definition}
\newtheorem{lemma}{Lemma}
\newtheorem{proposition}{Proposition}
\newtheorem{corollary}{Corollary}
\newtheorem{example}{Example}
\newtheorem{conjecture}{Conjecture}
\newtheorem{algorithm}{Algorithm}
\newtheorem{exercise}{Exercise}
\newtheorem{xample}{Example}
\newtheorem{remarkk}{Remark}
\newcommand{\bethm}{\begin{theorem}}
\newcommand{\ethm}{\end{theorem}}
\newcommand{\be}{\begin{equation}}
\newcommand{\ee}{\end{equation}}
\newcommand{\bea}{\begin{eqnarray}}
\newcommand{\eea}{\end{eqnarray}}
\newcommand{\beq}[1]{\begin{equation}\label{#1}}
\newcommand{\eeq}{\end{equation}}
\newcommand{\beqn}[1]{\begin{eqnarray}\label{#1}}
\newcommand{\eeqn}{\end{eqnarray}}
\newcommand{\beaa}{\begin{eqnarray*}}
\newcommand{\eeaa}{\end{eqnarray*}}
\newcommand{\req}[1]{(\ref{#1})}

\newcommand{\lip}{\langle}
\newcommand{\rip}{\rangle}
\newcommand{\uu}{\underline}
\newcommand{\oo}{\overline}
\newcommand{\La}{\Lambda}
\newcommand{\la}{\lambda}
\newcommand{\eps}{\varepsilon}
\newcommand{\om}{\omega}
\newcommand{\Om}{\Omega}
\newcommand{\ga}{\gamma}
\newcommand{\ka}{\kappa}
\newcommand{\rrr}{{\Bigr)}}
\newcommand{\qqq}{{\Bigl\|}}

\newcommand{\dint}{\displaystyle\int}
\newcommand{\dsum}{\displaystyle\sum}
\newcommand{\dfr}{\displaystyle\frac}
\newcommand{\bige}{\mbox{\Large\it e}}
\newcommand{\integers}{{\Bbb Z}}
\newcommand{\rationals}{{\Bbb Q}}
\newcommand{\reals}{{\rm I\!R}}
\newcommand{\realsd}{\reals^d}
\newcommand{\realsn}{\reals^n}
\newcommand{\NN}{{\rm I\!N}}
\newcommand{\DD}{{\rm I\!D}}
\newcommand{\degree}{{\scriptscriptstyle \circ }}
\newcommand{\dfn}{\stackrel{\triangle}{=}}
\def\complex{\mathop{\raise .45ex\hbox{${\bf\scriptstyle{|}}$}
     \kern -0.40em {\rm \textstyle{C}}}\nolimits}
\def\hilbert{\mathop{\raise .21ex\hbox{$\bigcirc$}}\kern -1.005em {\rm\textstyle{H}}} 
\newcommand{\RAISE}{{\:\raisebox{.6ex}{$\scriptstyle{>}$}\raisebox{-.3ex}
           {$\scriptstyle{\!\!\!\!\!<}\:$}}} 

\newcommand{\hh}{{\:\raisebox{1.8ex}{$\scriptstyle{\degree}$}\raisebox{.0ex}
           {$\textstyle{\!\!\!\! H}$}}}

\newcommand{\OO}{\won}
\newcommand{\calA}{{\mathcal A}}
\newcommand{\calB}{{\mathcal B}}
\newcommand{\calC}{{\cal C}}
\newcommand{\calD}{{\cal D}}
\newcommand{\calE}{{\cal E}}
\newcommand{\calF}{{\mathcal F}}
\newcommand{\calG}{{\cal G}}
\newcommand{\calH}{{\cal H}}
\newcommand{\calK}{{\cal K}}
\newcommand{\calL}{{\mathcal L}}
\newcommand{\calM}{{\mathcal M}}
\newcommand{\calO}{{\cal O}}
\newcommand{\calP}{{\cal P}}
\newcommand{\calT}{{\mathcal T}} 
\newcommand{\calU}{{\mathcal U}}
\newcommand{\calX}{{\cal X}}
\newcommand{\calY}{{\mathcal Y}}
\newcommand{\calZ}{{\mathcal Z}}
\newcommand{\calXX}{{\cal X\mbox{\raisebox{.3ex}{$\!\!\!\!\!-$}}}}
\newcommand{\calXXX}{{\cal X\!\!\!\!\!-}}
\newcommand{\gi}{{\raisebox{.0ex}{$\scriptscriptstyle{\cal X}$}
\raisebox{.1ex} {$\scriptstyle{\!\!\!\!-}\:$}}}
\newcommand{\intsim}{\int_0^1\!\!\!\!\!\!\!\!\!\sim}
\newcommand{\intsimt}{\int_0^t\!\!\!\!\!\!\!\!\!\sim}
\newcommand{\pp}{{\partial}}
\newcommand{\al}{{\alpha}}
\newcommand{\sB}{{\cal B}}
\newcommand{\sL}{{\cal L}}
\newcommand{\sF}{{\cal F}}
\newcommand{\sE}{{\cal E}}
\newcommand{\sX}{{\cal X}}
\newcommand{\R}{{\rm I\!R}}
\renewcommand{\L}{{\rm I\!L}}
\newcommand{\vp}{\varphi}
\newcommand{\N}{{\rm I\!N}}
\def\ooo{\lip}
\def\ccc{\rip}
\newcommand{\ot}{\hat\otimes}
\newcommand{\rP}{{\Bbb P}}
\newcommand{\bfcdot}{{\mbox{\boldmath$\cdot$}}}

\renewcommand{\varrho}{{\ell}}
\newcommand{\dett}{{\textstyle{\det_2}}}
\newcommand{\sign}{{\mbox{\rm sign}}}
\newcommand{\TE}{{\rm TE}}
\newcommand{\TA}{{\rm TA}}
\newcommand{\E}{{\rm E\,}}
\newcommand{\won}{{\mbox{\bf 1}}}
\newcommand{\Lebn}{{\rm Leb}_n}
\newcommand{\Prob}{{\rm Prob\,}}
\newcommand{\sinc}{{\rm sinc\,}}
\newcommand{\ctg}{{\rm ctg\,}}
\newcommand{\loc}{{\rm loc}}
\newcommand{\trace}{{\,\,\rm trace\,\,}}
\newcommand{\Dom}{{\rm Dom}}
\newcommand{\ifff}{\mbox{\ if and only if\ }}
\newcommand{\nproof}{\noindent {\bf Proof:\ }}
\newcommand{\remark}{\noindent {\bf Remark:\ }}
\newcommand{\remarks}{\noindent {\bf Remarks:\ }}
\newcommand{\note}{\noindent {\bf Note:\ }}

\newcommand{\boldx}{{\bf x}}
\newcommand{\boldX}{{\bf X}}
\newcommand{\boldy}{{\bf y}}
\newcommand{\boldR}{{\bf R}}
\newcommand{\uux}{\uu{x}}
\newcommand{\uuY}{\uu{Y}}

\newcommand{\limn}{\lim_{n \rightarrow \infty}}
\newcommand{\limN}{\lim_{N \rightarrow \infty}}
\newcommand{\limr}{\lim_{r \rightarrow \infty}}
\newcommand{\limd}{\lim_{\delta \rightarrow \infty}}
\newcommand{\limM}{\lim_{M \rightarrow \infty}}
\newcommand{\limsupn}{\limsup_{n \rightarrow \infty}}

\newcommand{\ra}{ \rightarrow }

\newcommand{\ARROW}[1]
  {\begin{array}[t]{c}  \longrightarrow \\[-0.2cm] \textstyle{#1} \end{array} }

\newcommand{\AR}
 {\begin{array}[t]{c}
  \longrightarrow \\[-0.3cm]
  \scriptstyle {n\rightarrow \infty}
  \end{array}}

\newcommand{\pile}[2]
  {\left( \begin{array}{c}  {#1}\\[-0.2cm] {#2} \end{array} \right) }

\newcommand{\floor}[1]{\left\lfloor #1 \right\rfloor}

\newcommand{\mmbox}[1]{\mbox{\scriptsize{#1}}}

\newcommand{\ffrac}[2]
  {\left( \frac{#1}{#2} \right)}

\newcommand{\one}{\frac{1}{n}\:}
\newcommand{\half}{\frac{1}{2}\:}

\def\le{\leq}
\def\ge{\geq}
\def\lt{<}
\def\gt{>}

\def\squarebox#1{\hbox to #1{\hfill\vbox to #1{\vfill}}}
\newcommand{\nqed}{\hspace*{\fill}
          \vbox{\hrule\hbox{\vrule\squarebox{.667em}\vrule}\hrule}\bigskip}

\title{Strong solutions of SDE with rough coefficients}

\author{Ali  S\"uleyman  \"Ust\"unel}

\begin{abstract}
  We give a proof of the strong existence and the regularity  of
  stochastic differential equations  driven by a Brownian motion and a measurable, Markovian drift without no regularity hypothesis except that the
  Girsanov exponential associated is in some $L^{1+\eps}(\mu)$ for
  some fixed $\eps>0$ by using the techniques which are totally novel
  originating from the abstract Wiener  space theory of Leonard Gross.
\noindent

\end{abstract}
\maketitle
\section{\bf{Introduction}}
\noindent
Let $W=C([0,1],\R^d)$ be the classical Wiener space, $H=\{h\in W:
h(t)=\int_0^t\dot{h}(s)ds,\dot{h}\in L^2([0,1],\R^d)\}$ is the
Cameron-Martin space and $\mu$ the Wiener measure associated. We
denote also by $(W_t,t\in [0,1])$ the canonical Wiener process
where $W_t$ is the evaluation function at $t$, i.e., for $w\in W$, $W_t(w)=w(t)$. Let $b:[0,1]\times \R^d\to\R^d$, be a measurable function, we assume the following:
\begin{itemize}
\item
  \begin{equation}
    \label{H1}
  E[\rho(-\delta \tilde{b})]=1\,,
  \end{equation}
  where
  $$
  \rho(-\delta\tilde{b})=\exp\left(-\int_0^1b(s,W_s).dW_s-\half\int_0^1|b(s,W_s)|^2ds\right)=e^{-f},
  $$
\item and  that there exists some constant $c_b$ such that,
  \begin{equation}
    \label{H2}
  E[\rho(-\delta\tilde{b})^{(1+c_b)}]=E[e^{-(1+c_b)f}]<\infty\,.
  \end{equation}
\end{itemize}

In this work we shall prove the strong existence and uniqueness of the SDE
$$
dX_t=b(t,X_t)dt+dW_t,\,X_0=x\in \R^d,
$$
without any further regularity hypothesis about the drift $b$, we
shall even show that the solution is real  $H$-analytic for $t<1$ by
transforming the problem to an abstract Wiener space, namely
$(\Om,H_2,P)$ where $\Om=C([0,1],W),\,H_2=H([0,1],H)\simeq
H\otimes_2H$ where $\otimes_2$ denotes the Hilbert-Schmidt tensor
product and $(B_t,t\in [0,1])$ is the evaluation map on $\Om$ which
happens to be a cylindrical Brownian motion under the Radonification
of the cylindrical Gaussian measure on $H_2$ as a Radon measure on
$(\Om,\calB(\Om))$, cf. \cite{L.G.,Kuo}. This  is  a totally different method from those employed in litterature (cf. \cite{Dav,F-F,K-R,M-P,R-Z}).

\section{Construction at the upper floor}
We mean by this title the existence an abstract Wiener space \`a la Leonard Gross: Let $H_2$ be the Hilbert space defined by
$$
H_2=\left\{k\in C([0,1],H): k(t)=\int_0^t\dot{k}(s)ds,\,\dot{k}(s)\in H\right\}\simeq H\otimes_2 H\,,
$$
with the norm $|k|_2^2=\int_0^1|\dot{k}(s)|_H^2ds$; note that this space is isomorphic to the Hilbert-Schmidt
tensor product of $H$ by itself. Let $\Om$ be $C([0,1],W)$ and $P$ be the Gaussian measure concentrated in $\Om$ with Cameron-Martin
space $H_2$ (cf. L. Gross and H.H. Kuo), denote by $B$ the coordinate process on $\Om$
under $P$ which is regarded as a $W$-valued Brownian motion, in fact
if we denote by $(B_{s,t}(\om),s,t\in [0,1])$ the family
$(B_s(\om)(t),s,t\in [0,1])$ we obtain the Brownian sheet. In
particular the law of $B_s$ under $P$, denoted as $\mu_s$ is the Gauss
measure on $W$ whose the characteristic function is
$$
\hat{\mu}_s(h)=\exp-\frac{s}{2}|h|_H^2\,.
$$
Let us note that $B$ is an additive process with semigroup $Q_t$ defined by $Q_tG(w)=E_P[G(w+B_t)]$ with $G:\Om\to\R_+$ for
measurable function $G$. From the It\^o formula at the upper floor, we can represent the martingale $(E_P[e^{-f\circ B_1}|\calB_t],t\geq 0)$ as a
stochastic integral, at the same time, as $B$ is an additive process we have
$$
E_P[e^{-f\circ B_1}|\calB_t]=Q_{1-t}(e^{-f})(B_t)\,,
$$
where $Q_t$ is the heat kernel associated to $B$, i.e.,
$$
Q_tF(w)=\int_W F(w+y)\mu_t(dy)=E_P[G(B_t+w)]\,,
$$
where $\mu_t$ is the Gauss measure on $W$ with variance $tI_H$ as
explained above. Hence from the It\^o formula, we have
\beaa
Q_{1-t}(e^{-f})(B_t)&=&1+\int_0^t(\nabla (Q_{1-s}(e^{-f})(B_s),dB_s)\\
&=&1+\int_0^t Q_{1-s}(e^{-f})(B_s)(\nabla\log Q_{1-s}(e^{-f})(B_s),dB_s)\,,
\eeaa
where the stochastic integral is the well-known cylindrical stochastic
integral.  Hence $Q_{1-t}(e^{-f})(B_t)$ can be represented as
$$
Q_{1-t}(e^{-f})(B_t)=\exp\left(-\int_0^t(\dot{v}_s(B_s),dB_s)-\half\int_0^t|\dot{v}_s(B_s)|_H^2ds\right)
$$
where $\dot{v}_t(B_t)=-\nabla \log Q_{1-t}e^{-f}(B_t)$ and  the stochastic integral is the classical cylindrical  stochastic integral, regarding $B$ as a cylindrical Brownian motion indexed by $H$,
well studied (cf. \cite{Kuo} for instance) during these last 50 years.

Before proceeding further, we shall give a modification of the Paley-Wiener integral which will be well-defined for any Gaussian random variable
$B_s$, spontaneously for any $s\geq 0$ as suggested in \cite{K-S}: if $h\in W^\star, w\in W$, we define $I(h)(w)$ as $<h,w>$ where $<\cdot,\cdot>$ is the duality bracket
corresponding to the dual pair $(W^\star,W)$, if not, let $(z_n,n\geq 1)\subset W^\star$ be a sequence converging to $h$ (in $H$), define
$$
D(h)=\{w\in W: \lim_n<z_n,w>\text{ exists}\}\,,
$$
then define $I(h)(w)$ as the limit of $(<z_n,w>,n\geq 1)$ if $w\in D(h)$ and as zero if not. Note that $D(h)+H\subset D(h)$, hence $D(h)^c$ is a slim
set.
After this construction, we have the following observation: if $h\in H$, then one has
\bea
  \label{chv}
  Q_{1-t}(e^{-f})(w+h)&=&\int_W e^{-f(w+h+y)}\mu_{1-t}(dy)\\
  &=&\int_W e^{-f(w+y)}\exp\frac{1}{1-t}(I(h)(y)-\half |h|_H^2)\mu_{1-t}(dy)\nonumber
  \eea
  From the expression (\ref{chv}) it follows that the function $w\to Q_{1-t}(e^{-f})(w)$ is real $H$-analytic, i.e., $\mu_t$-almost all $w\in W$, the map
  $h\to Q_{1-t}(e^{-f})(w+h)$ is real analytic on $H$ for ant $t\in [0,1)$.
    Here is our important theorem:
   \bethm
     Let $\tau<1$ be a fixed number and define an adapted perturbation of identity $V^\tau:\Om\to \Om$ on the upper floor as
     $$
     V^\tau_t(B)=B_t(B)+\int_0^{\tau\wedge t} \dot{v}_s(B_s(B))ds\,.
     $$
     Then $V^\tau=I_B+v^\tau$ is $P$-a.s. invertible.
     \ethm
     \nproof From \ref{chv}, it follows that $h\to \dot{v}_s(w+h)$ is infinitely differentiable $\mu_s$-a.s. on $H$. Moreover, from martingale property, we have
     $$
     \inf_{t<1}Q_{1-t}(e^{-f})(B_t)>0
  $$
     $P$-a.s.
     Hence by $H$-regularity and $t$-regularity of $\dot{v}$ on $[0,\tau]$, we deduce
     \begin{equation}
       \label{cond_1}
       \sup_{t\leq \tau}\sup_{\|K\|_2\leq M}|\dot{v}_t(B_t(\om)+K_t)|_H\leq c_{1,M}(\tau,\om)<\infty\,,
       \end{equation}
     \begin{equation}
       \label{cond_2}
       \sup_{t\leq \tau}\sup_{\|K\|_2\leq N}\|\nabla\dot{v}_t(B_t(\om)+K_t)\|_{H_2}\leq c_{2,N}(\tau,\om)<\infty\,,
     \end{equation}
     for any $M,N>0$, where $c_{1,M}(\tau,\om)$ and
     $c_{2,N}(\tau,\om)$ are $P$-almost surely  finite random variables. From these two estimates follows direcly the $H-C^1$-property of $v^\tau$ on
     the abstract Wiener space $(\Om,H_2,P)$.

     Recall that $(\dot{v}_t(B_t(B)),\,t\in [01,1])$ is adapted to the filtration of the canonical Brownian motion $(B_t)$ at the upper floor and
     $$
     E_P[\rho(-\delta_Bv^\tau)]=E_P[\rho_\tau(-\delta_Bv)]=E[e^{-f}]=1.
     $$
     By $H$-regularity, it comes from the change of variables formula that,
     for any $G\in C_b(\Om)$, we have by the change of variables formula on an abstract Wiener space, (cf.\cite{BOOK}, Chapter 3)
     \begin{equation}
       \label{chg-2}
     E_P[G\circ V^\tau \La]=E_P[G\,N(\cdot,V^\tau)]=E_P[G]\,,
     \end{equation}
     where 
     $\La=\dett(I_{H_2}+\nabla v^\tau)\rho(-\delta v^\tau)=\rho(-\delta
     v^\tau)$ and  $\dett(I_{H_2}+\nabla v^\tau)$ denotes the modified Carlman-Fredholm
     determinant of $I_{H_2}+\nabla v^\tau$ which is equal to one because
     of the adaptedness of $(t,B)\to\dot{v}^\tau(t,B_t(B))$ and finally $N(\om,V^\tau)$ denotes the multiplicity
     of $V^\tau$, at $\om\in \Om$, i.e., the cardinality of $(V^\tau)^{-1}\{\om\}$, which happens to be equal to one $P$-almost everywhere due to the
     relation (\ref{chg-2}) and this suffices to conclude  the almost sure invertibility of $V^\tau$.
     \nqed
     \bethm
     \label{upper_soln}
     The following stochastic differential equation
     \begin{equation}
     \label{SDE_1}
     dU_t=-\dot{v}(t,U_t)dt+dB_t,\,U_0=0,\, t\in[0,1]
     \end{equation}
     has unique strong solution on $(\Om,\calB(\Om),P)$.
     \ethm
     \nproof
     Let $U^\tau$ be the both-sided inverse of $V^\tau$, it satisfies
     \beaa
     dU^\tau&=&-\dot{v}^\tau(t,U^\tau_t)dt+dB_t\\
     &=&-\dot{v}(t,U^\tau_t)\,1_{[0,\tau]}(t)dt+dB_t\\
     \eeaa
     If $\kappa<\tau$, then it is easy to see that $U^\tau_{t\wedge\kappa}=U^\kappa_{t\wedge \kappa}$, hence the processes $(U^\tau,\tau<1)$
     can be pieced out by defining $U$ as $U_t=U^\tau_t$ if $t\leq
     \tau$ for some $\tau<1$.  Let $(\Pi_n,n\geq 1)$ be a sequence of
     finite, increasing subsets of $[0,1)$, as
     $$
     \frac{dU^\tau(P)}{dP}=\rho(-\delta_Bv^\tau)=\rho_\tau(-\delta_Bv)\,.
     $$
    we get
       \beaa
       E_P[\sup_{t<1}\|U_t\|_W]&=&\lim_nE_P[\sup_{t\in \Pi_n}\|U_t\|_W]\\
       &=&\lim_nE_P[\sup_{t\in \Pi_n}\|B_t\|_W\rho(-\delta_B(v)]\\
       &=&E_P[\sup_{t<1}\|B_t\|_W\rho(-\delta_B(v)]\\
       &=&E_P[\sup_{t\leq 1}\|B_t\|_W e^{-f\circ B_1}]<\infty\,.
     \eeaa
     Moreover
     $$
     \lim_{s,t\to 1}E_P[\|U_t-U_s\|_W]=\lim_{s,t\to 1}E_P[\|B_t-B_s\|_W e^{-f\circ B_1}]=0\,,
     $$
     hence $\lim_{t\to 1} U_t$ does exist, we denote it naturally by $U_1$. Moreover
     $$
     U_1\circ V=\lim_{t\to 1}U_t\circ V=\lim_{t\to 1}B_t=B_1\,,
     $$
     hence $U_1$ belongs to the inverting family $(U_t)$ of $V$, in other words $U=(U_t,t\in [0,1])$ is the almost sure inverse of $V$, moreover
     by the martingale convergence theorem, we get
     $$
     \frac{dU(P)}{dP}=\rho(-\delta_Bv)=e^{-f\circ B_1}\,.
     $$
     It is trivial now to see that $U$ is the strong solution of the SDE (\ref{SDE_1}).
     \nqed

     \section{Variational identification}
     \noindent
     Let us denote by $L_a^2(P,H_2)$ the obvious completion of the  set of maps $\xi:\Om\to H_2$ such that $\xi(\om)(t)=\int_0^t\dot{\xi}(s,\om)ds$ with
     $\dot{\xi}(s):\Om\to H$ is $\calB_s$ measurable for almost all $s\in [0,1]$ such that
     $$
     \|\xi\|_{2,a}^2= E_P\|\xi\|_{H_2}^2=\int_0^1|\dot{\xi}(s)|_H^2ds<\infty\,,
     $$
     and the completion mentioned above is to be understood under this norm of the equivalnce classes. In particular the loer index ``a'' refers to the
     causal character of $\dot{\xi}$.

     We know that the strong existence of the SDE (\ref{SDE_1}) is equivalent to the existence of the unique minimizing element of the
     following variational problem (cf.,
     \cite{ASU-2,ASU-4})\footnote{In fact in these references the
       classical Wiener space has been treated but the results extend
       easily to our abstract Wiener space $(\Om,H_2,P)$ with the same
     proof due to stochastic integral representation of the
     $(\calB_t,t\in [0,1])$-martingales.}:
     \begin{equation}
       \label{upper}
       \inf( K(\xi): \xi\in L^2_a(P,H_2))=\inf \left[E_P\left[\half \|\xi\|_{H_2}^2 +f\circ B_1\circ (I_B+\xi)\right]:\xi\in L^2_a(P,H_2)\right]=0\,.
     \end{equation}
     Let us calculate the functional $K$ in more detail: recall that we have used the notation $\tilde{b}(w)=\int_0^\cdot b(s,W_s(w))ds$ and  also
     $f(w)=\delta \tilde{b}(w)+\half|\tilde{b}(w)|_H^2$. We then have
     $$
     K(\xi)=E_P\left[\half\|\xi\|_{H_2}^2+\left(\delta\tilde{b}+\half|\tilde{b}|^2\right)\circ (\xi_1(B)+B_1(B))\right]
     $$
     where $B$ is the generic point of $\Om$, which is $W$-valued Brownian  path under  $P$ used for pedagogical convenience. The second term at the
     right can be written as:
     \beaa
     E_P[\delta\tilde{b}(\xi_1(B)+B_1)]&=&E_P\left[\left(\int_0^1 b(s,W_s)dW_s\right)\circ(\xi_1(B)+B_1)\right]\\
     &=& E_P\left[\int_0^1b(s,W_s(\xi_1(B)+B_1)).(d(W_s(B_1(B)+B_1)\right]\\
     &=&E_P\left[\int_0^1b(s,\xi_{1,s}(B)+B_1(s)).(\dot{\xi}_{1,s}ds+B_1(ds))\right]\\
     &=&E_P\left[\int_0^1b(s,\xi_{1,s}(B)+B_1(s)).\,\dot{\xi}_{1,s}ds\right]\\
     &=&E_P\left[\int_0^1\int_0^1b(s,\xi_{1,s}(B)+B_1(s)).\,\ddot{\xi}_{t,s}dt ds\right]\,
     \eeaa
     where the term with the stochastic integrator $B_1(ds)$ does not contribute because of the expectation operation since $B_{1,s_{i+1}}-B_{1,s_i}$ is
     orthogonal to $b(s_i,\xi_{1,s_i}(B)+B_{1,s_i})$ in $L^2(P)$ for smooth $b$ and this property perpetuates under $L^2$-limit. Hence we obtain
     \beaa
     K(\xi)&=& E_P\left[\half \|\xi\|_{H_2}^2+\int_{[0,1]^2}b(s,\xi_{1,s}(B)+B_1(s))\cdot \ddot{\xi}_{t,s}(B)dsdt\right.\\
     &&\left.+\int_0^1|b(s,\xi_{1,s}(B)+B_1(s))|^2ds\right]\\
     &=&\half
     E_P\left[\|\xi\|_{H_2}^2+\int_0^1|\dot{\xi}_{1,s}+b(s,\xi_{1,s}(B)+B_1(s))|^2ds\right]-\half
     E_P\int_0^1|\dot{\xi}_{1,s}|^2ds\\
     &=&E_P\left[\int_0^1|\dot{\xi}_{1,s}+b(s,\xi_{1,s}(B)+B_1(s))|^2ds\right]\,.
       \eeaa
     If $\xi$ is the unique minimizer of $K$, since $E[e^{-f}]=1$,
     letting, $X_t=\xi_{1,t}+B_1(t)$ we realize that $X$ satisfies the
     SDE
     $$
     dX_t=-b(t,X_t)dt+dB_{1,t}
     $$
     almost surely.
     \begin{theorem}
       \label{strong}
       The process $X=(X_t,t\in [0,1])$ is adapted to the filtration
       of  $B_1$.
       \end{theorem}
\nproof
     We claim that this solution is the strong one.  To
     prove this we shall almost repeat the same strategy that we have
     employed upto now, the only difference being the
     modification-regularization of $e^{-f}$. Namely define $f_n$ by
     \begin{equation}
       \label{BN}
     e^{-f_n}=P_{1/n}E[e^{-f}|V_n]=\rho(-\delta b_n)
     \end{equation}
     where $P_{1/n}$ is the Ornstein-Uhlenbeck semigroup on $W$, $V_n$
     is the sigma algebra generated by
     $$
     \{w\to <z_i,w>;\,i=1,\ldots,n\}
     $$
     where $z_i\in W^\star$  and the images
     of $z_i$ in $H$ are orthonormal and $b_n$ is defined via Ito
     representation theorem. By  the hypothesis
     $E[e^{-(1+c_b)f}]<\infty$, $f_n$ is an $H-C^\infty$-map (cf. \cite{KH})
     consequently, we can represent the martingale
     $(E_P[e^{-f_n(B_1)}|\calB_t]=Q_{1-t}e^{-f_n}(B_t),t\in [0,1])$ as
       an exponential martingale
       $$
       Q_{1-t}e^{-f_n}(B_t)=\rho_t(-\delta_Bv_n)=\exp\left(-\int_0^t(\dot{v}_n(s),dB_s)-\int_0^t|\dot{v}_n(s)|_H^2ds\right)
       $$
       thanks to the Ito formula for the cylindrical Brownian motion
       $(B_t,t\in [0,1])$, where $w\to \dot{v}_n(t,w)=\nabla \log
       Q_{1-t}(e^{-f_n})(w)$ is an $H-C^\infty$-map. Therefore the
       perturbation of identity defined at the upper floor by
       $$
       B\to V_n(B)(t)=B(t)+\int_0^t \dot{v}_{n}(s,B_s(B))ds
       $$
     is an $H-C^\infty$-map on the Cameron-Martin space of the upper
     floor, namely $H_2$ (or the upper floor abstract Wiener space
     $(\Om,H_2,P)$). Using the same argument as we have done for
     $V^\tau$, it follows that $V$ is $P$-a.s. invertible with inverse $U_n$
     satifying the strong  SDE
     $$
     dU_n(t)=-\dot{v}_n(t,U_n(t))dt+dB_t\,,
     $$
     where $(t,\om)\to u_n(t,\om)=-\int_0^t \dot{v}_n(s,U_n(s,\om))ds$
     is the unique solution of the minimization problem for the
     functional on $L^2_a(P,H_2)$ defined by
     $$
     \eta\to E_P\left[\half\|\eta\|_2^2+f_n\circ B_1\circ(I_\Om+\eta)\right]=K_n(\eta)\,.
     $$
     This implies, as we have already seen above, that $U_n$ satisfies
     \begin{equation}
       \label{var_eq}
     dU_n^1(t)=-b_n(t,U_n^1(t))dt +dB^1_t\,.
     \end{equation}
As the vector field $b_n$ defined by the relation \ref{BN} is smooth,
the SDE
$$
dX_t=-b_n(t,X_t)dt+dW_t
$$
has unique strong solution in the lower floor Wiener space, i.e.,
$(t,w)\to X_t(w)=F(t,w)$,  $F$ being adapted to $(\calF_t)$ which is
the filtration of $(W_t,t\in [0,1])$. This
implies that the weak solution of the SDE \ref{var_eq} is a
strong solution, in other words $(U_n(t))$ is adapted to the filtration
of $B^1=(B^1_t, t\in [0,1])$. Recall now that we have
$$
\dot{v}_n(t,B_t)=-\nabla \log Q_{1-t}e^{-f_n}(B_t)\to \dot{v}(t,B_t)
$$
$P$-a.s., as $n\to\infty$. Let now $g$ and $H$ be two nice, bounded, continuous
functions on $W$, assume further that $g$ is $\calF_t$-measurable. We have
\beaa
\lim_n E_P[g\circ U_n^1\,H(B_1)]&=&\lim_n E_P[g\circ
B_1\,H(V^n_1)e^{-f_n\circ B_1}]\\
&=&E_P[g\circ B_1 H(V_1)e^{-f\circ B_1}]\\
&=&E_P[g\circ U^1H(B_1)]\,,
\eeaa
consequently the sequence $(g\circ U_n,n\geq 1)$ converges weakly in
any $L^p(P)$ to $g\circ U^1$, hence a subsequence of its convex
combinations converges strongly, therefore $g\circ U^1$ is measurable
w.r.t. sigma algebra generated by $(B^1_s,s\in [0,t])$ for any $g$ as above, hence
$U^1$ itself also is adapted to the filtration of $B^1$, hence $U^1$
is a strong solution.
\nqed

\vspace{2cm}

{\footnotesize{\bf{
\noindent
A.S. \"Ust\"unel, 70, rue Bobillot, 75013 Paris, France\\
or\\
M\"unir Nur. Sel\c{c}uk Cadd. No.16/18, Fenerbah\c{c}e, Istanbul, Turkey,\\
asustunel@gmail.com
 
}}

\end{document}